\documentclass[a4paper,11pt]{article}

\usepackage[T1]{fontenc}
\usepackage{lmodern}
\usepackage{amsmath,amssymb,amsfonts}
\usepackage{geometry}
\usepackage{array,longtable}
\usepackage[hidelinks]{hyperref}
\usepackage{microtype}

\allowdisplaybreaks
\newcommand{\Z}{\mathbb Z}
\newcommand{\Bcal}{\mathcal B}
\newcommand{\Dcal}{\mathcal D}
\newcommand{\sr}{\mathrm{sr}}
\newcommand{\ar}{\mathrm{ar}}

\hypersetup{
  pdftitle={On Borwein polynomials with multiple roots outside the unit disk},
  pdfauthor={Souad El Otmani; Musbahu Idris; Jean-Marc Sac-Epee},
  pdfkeywords={Borwein polynomials; multiple roots; Lehmer's polynomial; integer linear programming}
}

\title{On Borwein polynomials with multiple roots outside the unit disk}

\author{Souad El Otmani\thanks{Institut \'Elie Cartan de Lorraine, Universit\'e de Lorraine, France.\newline
Email: \texttt{souad.elotmani@univ-lorraine.fr}}
\and Musbahu Idris\thanks{Institut \'Elie Cartan de Lorraine, Universit\'e de Lorraine, France.\newline
Email: \texttt{musbahu.idris@univ-lorraine.fr}}
\and Jean-Marc Sac-\'Ep\'ee\thanks{Institut \'Elie Cartan de Lorraine, Universit\'e de Lorraine, France.\newline
Email: \texttt{jean-marc.sac-epee@univ-lorraine.fr}}}
\date{}

\begin{document}

\maketitle

\begin{abstract}
A Borwein polynomial has coefficients in \(\{-1,0,1\}\) and nonzero
constant and leading coefficients. We seek Borwein polynomials having a multiple root outside the unit
disk. To this end, we consider the powers \(L^2\) and \(L^3\) of
Lehmer's polynomial \(L\), which has a real root \(\tau>1\).
Since \(L\) has exceptionally small Mahler measure, the theorem of
Bombieri and Vaaler~\cite{BombieriVaaler1987} guarantees the existence
of Borwein multiples of these powers. Such multiples have \(\tau\) as
a root of multiplicity at least \(2\) and \(3\), respectively. Mossinghoff~\cite[Algorithm~3,
p.~323]{Mossinghoff2003} constructed examples that, under the substitution
\(x\mapsto -x\), yield such multiples.
The main contribution of this paper is to determine the least possible
degrees of Borwein multiples of \(L^2\) and \(L^3\), with and without
symmetry restrictions on their coefficients, and to describe their degree spectra. For \(L^2\), the minimum degree is \(34\) without symmetry
restrictions, \(38\) among self-reciprocal multiples, and \(41\) among
antireciprocal multiples. For \(L^3\), the corresponding minimum degrees
are \(115,115,117\). We determine the complete degree spectra for \(L^2\) and, for \(L^3\),
classify the degrees from \(115\) to \(150\) for which self-reciprocal or
antireciprocal Borwein multiples exist, giving an explicit example for each
admissible degree in each symmetry class. Further constructions
yield Borwein multiples of \(L^3\) in every degree \(D\geq215\), with
both symmetry types for every \(D\geq231\).
\end{abstract}



\section{Introduction}

Following Drungilas, Jankauskas, and
\v{S}iurys~\cite{DrungilasJankauskasSiurys2018}, a \emph{Borwein
polynomial} of degree \(D\geq0\) is a polynomial
\[
 Q(x)=\sum_{j=0}^{D}q_jx^j\in\Z[x],
 \qquad q_j\in\{-1,0,1\},
 \qquad q_0q_D\ne0.
\]
\begin{samepage}
Let \(P\in\mathbb C[x]\) be a nonzero polynomial, and write
\(P(x)=a_d\displaystyle\prod_{j=1}^{d}(x-\alpha_j)\).
Its Mahler measure is
\[
 M(P)=|a_d|\prod_{j=1}^{d}\max\{1,|\alpha_j|\}.
\]
\end{samepage}
Lehmer's problem~\cite{Lehmer1933,Smyth2008} asks whether there exists a
constant \(c>1\) such that \(M(P)\geq c\) for every \(P\in\Z[x]\)
with \(M(P)>1\).

For the question studied here, the relevance of small Mahler measure comes
from Corollary~2 of Bombieri and Vaaler~\cite{BombieriVaaler1987}, which
implies that every monic polynomial \(R\in\Z[x]\) with \(M(R)<2\)
has a nonzero multiple in \(\Z[x]\) of height \(1\), that is, with coefficients in
\(\{-1,0,1\}\).  If \(R(0)\ne0\), dividing such a multiple by the largest
power of \(x\) that divides it yields a Borwein multiple of \(R\).

Let \(P\) be the minimal polynomial of an algebraic integer
\(\alpha\) with \(|\alpha|>1\), and let \(m\geq1\).  The polynomial \(P^m\)
is monic and has nonzero constant term.  Moreover, by the multiplicativity of
Mahler measure,
\[
 M(P^m)=M(P)^m.
\]
Since \(P\) is the minimal polynomial of \(\alpha\), a polynomial
\(Q\in\Z[x]\) has \(\alpha\) as a root of multiplicity at least \(m\)
if and only if \(P^m\) divides \(Q\). Consequently, by the result
of Bombieri and Vaaler recalled above, the condition \(M(P)^m<2\)
guarantees the existence of a Borwein polynomial having \(\alpha\) as
a root of multiplicity at least \(m\).

Earlier work on multiples of irreducible integer polynomials includes
that of Pathiaux~\cite{Pathiaux1972} and Mignotte~\cite{Mignotte1975}.
An example directly related to multiple roots was given by Drungilas,
Jankauskas, and \v{S}iurys~\cite[Example~9]{DrungilasJankauskasSiurys2018}:
for \(P(x)=x^3-x+1\), they construct the degree-\(8\) Borwein multiple
of \(P^2\)
\[
\begin{aligned}
 (x^2+x+1)P(x)^2
 &=(x^2+x+1)(x^3-x+1)^2\\
 &=x^8+x^7-x^6+x^4+x^3-x+1,
\end{aligned}
\]
which has a double real root outside the unit disk.

The condition \(M(P)^m<2\) suggests choosing \(P\) with small
Mahler measure when seeking Borwein polynomials with multiple roots
outside the unit disk. We therefore turn to Lehmer's polynomial
\(L\)~\cite{Lehmer1933}:
\begin{equation}
 L(x)=x^{10}+x^9-x^7-x^6-x^5-x^4-x^3+x+1.
\end{equation}
The polynomial \(L\) is
irreducible over \(\mathbb Q\).  Its roots consist of two real roots off the
unit circle, \(\tau>1\) and \(\tau^{-1}\), and eight roots on the unit circle;
see Smyth~\cite{Smyth2008}.  Hence
\[
 M(L)=\tau=1.176280818\ldots.
\]
This remains the smallest known Mahler measure greater than \(1\) of an
integer polynomial; see El-Serafy and McKee~\cite[p.~603]{ElSerafyMcKee2025}.
In particular,
\[
 M(L)^4=1.914450139\ldots<2.
\]
Thus, for \(m=2,3,4\), \(M(L)^m<2\), and \(L^m\) has a Borwein multiple.

Neither the least degree nor the set of possible degrees of these multiples
follows from Corollary~2 of Bombieri and Vaaler.  Since \(L\) itself is a
Borwein polynomial, the first two nontrivial cases are \(L^2\) and \(L^3\).

For a polynomial \(F\) of degree \(d\), set
\[
 F^*(x)=x^dF(1/x).
\]
A polynomial \(F\) is \emph{self-reciprocal} if \(F^*=F\), and
\emph{antireciprocal} if \(F^*=-F\).  Mossinghoff~\cite[p.~315]{Mossinghoff2003}
uses the term \emph{reciprocal} for both signs and considers
\[
 \ell(x)=x^{10}-x^9+x^7-x^6+x^5-x^4+x^3-x+1=L(-x).
\]
Both \(L\) and \(\ell\) are self-reciprocal.  He constructs a
self-reciprocal Borwein polynomial of degree \(38\) divisible by \(\ell^2\),
an antireciprocal Borwein polynomial of degree \(115\) divisible by
\(\ell^3\), and a self-reciprocal Borwein polynomial of degree \(117\)
also divisible by \(\ell^3\)
\cite[Algorithm~3, p.~323]{Mossinghoff2003}.

The substitution \(x\mapsto -x\) gives a degree-preserving bijection
between Borwein multiples of \(L^m\) and those of \(\ell^m\).  It preserves
self-reciprocity and antireciprocity in even degree and interchanges them in
odd degree.  Thus these examples
yield self-reciprocal Borwein multiples of \(L^2\) and \(L^3\) of degrees
\(38\) and \(115\), respectively, as well as an antireciprocal Borwein
multiple of \(L^3\) of degree \(117\).

Let \(\Bcal\) denote the set of Borwein polynomials.  For \(m\geq1\), define
\begin{align*}
 \Dcal_m(L)&=\{\deg Q:Q\in\Bcal,\ L^m\mid Q\},\\
 \Dcal_m^{\sr}(L)&=\{\deg Q:Q\in\Bcal,\ L^m\mid Q,\ Q^*=Q\},\\
 \Dcal_m^{\ar}(L)&=\{\deg Q:Q\in\Bcal,\ L^m\mid Q,\ Q^*=-Q\}.
\end{align*}

For \(m=2,3\), let the least degrees in the unrestricted and two
symmetry-restricted classes be denoted by
\[
 \delta_m(L)=\min\Dcal_m(L),\qquad
 \delta_m^{\sr}(L)=\min\Dcal_m^{\sr}(L),\qquad
 \delta_m^{\ar}(L)=\min\Dcal_m^{\ar}(L).
\]

Our results determine these six minimum degrees:
\[
\begin{aligned}
 \delta_2(L)&=34, &\delta_2^{\sr}(L)&=38, &\delta_2^{\ar}(L)&=41,\\
 \delta_3(L)&=115, &\delta_3^{\sr}(L)&=115, &\delta_3^{\ar}(L)&=117.
\end{aligned}
\]
We also determine the complete unrestricted, self-reciprocal, and
antireciprocal degree spectra for \(L^2\). For \(L^3\), we classify the
unrestricted degrees through \(118\) and the degrees in each symmetry
class from \(115\) to \(150\). Explicit constructions provide further
admissible degrees and show that every degree \(D\geq215\) occurs in the
unrestricted case, while every degree \(D\geq231\) occurs in both
symmetry classes. The full degree spectrum for \(L^3\) is not yet known.

For a prescribed degree \(D\) and \(m\in\{2,3\}\), we seek a Borwein
polynomial \(Q\) of degree \(D\) divisible by \(L^m\). We write
\(Q=L^mH\), with \(H\in\Z[x]\), and regard the coefficients of \(H\)
as unknown integers. Each coefficient of \(Q\) is then a linear form in
these unknowns. Requiring these forms to lie between \(-1\) and \(1\),
with nonzero constant and leading coefficients, gives an integer
feasibility problem for each prescribed degree.

We test successive degrees, starting at \(D=\deg(L^m)=10m\). When the
system has no solution, that degree is excluded; a solution yields a
Borwein multiple of that degree. The first degree admitting a solution
is therefore the least possible one. Testing higher degrees also reveals
which degrees are admissible, including possible gaps above the minimum.
We describe the systems, with and without symmetry restrictions, in the
next section.

\section{Integer-linear programming for Borwein multiples}
\label{sec:models}

\subsection{The feasibility system for multiples without symmetry restrictions}
\label{subsec:unrestricted-model}

Fix \(m\in\{2,3\}\), and let \(D\) be the prescribed degree of a Borwein polynomial \(Q\) divisible by \(L^m\). Since \(L\) has degree \(10\), necessarily \(D\geq10m\). Write
\[
 L(x)^m=\sum_{i=0}^{d}g_ix^i,
 \qquad d=10m,
 \qquad e=D-d.
\]
The integers \(g_i\) are known.  We seek integers \(h_0,\ldots,h_e\) such
that, for
\[
 H(x)=\sum_{j=0}^{e}h_jx^j,
\]
the product \(Q=L^mH\) is a Borwein polynomial of degree \(D\).  Set
\(h_j=0\) when \(j\notin\{0,\ldots,e\}\).  The coefficient of \(x^k\) in
\(Q\) is
\[
 \sum_{i=0}^{d}g_i h_{k-i}.
\]
Since this coefficient is an integer, constraining it to the interval
\([-1,1]\) forces it to belong to \(\{-1,0,1\}\).  Moreover,
\(g_0=g_d=1\), so the constant and leading coefficients of \(Q\) are \(h_0\)
and \(h_e\), respectively.  Replacing \(H\) by \(-H\) if necessary permits
the normalization \(h_0=1\).  We thus obtain the following system:
\begin{equation}\label{eq:unrestricted-system}
\begin{array}{ll}
\text{find}
 & (h_0,\ldots,h_e)\in\Z^{e+1},\\[2mm]
\text{subject to}
 & h_0=1,\qquad h_e\in\{-1,1\},\\
 & \displaystyle -1\leq\sum_{i=0}^{d}g_i h_{k-i}\leq1
  \qquad(1\leq k<D).
\end{array}
\end{equation}
For each of the two choices \(h_e=1\) and \(h_e=-1\),
\eqref{eq:unrestricted-system} is an integer linear programming (ILP)
feasibility problem.  We solve these two cases separately.  The system
has a solution if and only if \(L^m\) has a Borwein multiple of degree
\(D\); the multiple is then \(Q=L^mH\).

\subsection{Feasibility systems for self-reciprocal and antireciprocal multiples}
\label{subsec:reciprocal-model}

To treat both symmetry types, note that \(L^m\) is self-reciprocal.
For \(Q=L^mH\), we therefore have \(Q^*=Q\) if and only if \(H^*=H\),
and \(Q^*=-Q\) if and only if \(H^*=-H\).
In the self-reciprocal case, we add to
\eqref{eq:unrestricted-system} the equations
\[
 h_j=h_{e-j}
 \qquad(0\leq j\leq e).
\]
In particular, \(h_e=h_0=1\).  No parity assumption is needed: if \(e\) is
odd, the coefficients of \(H\) occur in symmetric pairs; if \(e\) is even,
\(h_{e/2}\) is the central coefficient.  The resulting problem has a solution
if and only if \(L^m\) has a self-reciprocal Borwein multiple of degree \(D\).

In the antireciprocal case, we instead add
\[
 h_j=-h_{e-j}
 \qquad(0\leq j\leq e)
\]
to \eqref{eq:unrestricted-system}.  Thus \(h_e=-h_0=-1\), and when \(e\) is
even the central coefficient is \(h_{e/2}=0\).  This system has a solution if
and only if \(L^m\) has an antireciprocal Borwein multiple of degree \(D\).

\subsection{Degree-by-degree computations and exact verification}

We used Gurobi~\cite{Gurobi} to determine the feasibility of
system~\eqref{eq:unrestricted-system} for each \(m\in\{2,3\}\) and
every degree in the finite ranges specified below. In the
symmetry-restricted ranges, we also determined the feasibility of
the corresponding self-reciprocal and antireciprocal systems.
Infeasibility excludes the degree in the corresponding class; a solution
gives a polynomial \(H\), and hence \(Q=L^mH\). For every solution found,
we checked in exact integer arithmetic that \(Q\) is a Borwein polynomial
of degree \(D\). In the symmetry-restricted cases, we also verified
that \(Q\) has the required symmetry.

\section{Borwein multiples of \texorpdfstring{$L^2$}{L squared}}

We fix \(m=2\) and solve the ILP feasibility problem
\eqref{eq:unrestricted-system} for each \(20\leq D\leq68\), and its
symmetry-restricted variants for each \(20\leq D\leq76\).
As we show in the following subsections, these cutoffs suffice: the degree-\(34\) multiple \(Q_{34}\) and the degree-\(38\)
self-reciprocal multiple \(Q_{38}\) yield, by disjoint-support shifts,
multiples without symmetry restrictions in every degree \(D\geq69\) and
multiples of both symmetry types in every degree \(D\geq77\), respectively.

\subsection{The least degree among multiples without symmetry restrictions}

The unrestricted feasibility problems have no solution in degrees
\(20\leq D\leq33\). At degree \(34\), the computation yields a polynomial
\(H_{14}\) whose coefficients satisfy
system~\eqref{eq:unrestricted-system}. The corresponding Borwein polynomial is
\[
 Q_{34}(x)=L(x)^2H_{14}(x),
\]
where
\begin{align*}
 H_{14}(x)={}&x^{14}-3x^{13}+4x^{12}-3x^{11}+3x^9-4x^8+3x^7\\
 &-2x^6+x^5-x^3+2x^2-2x+1.
\end{align*}
The full coefficient vector of \(Q_{34}\) is given in
Table~\ref{tab:l2-unrestricted-polynomials} of the appendix.
Direct multiplication shows that every coefficient of \(Q_{34}\) belongs to
\(\{-1,0,1\}\), while its constant and leading coefficients are nonzero.
It follows that
\begin{equation}
 \delta_2(L)=34.
\end{equation}
The polynomial \(Q_{34}\) is neither self-reciprocal nor antireciprocal.

\subsection{The least degrees among self-reciprocal and antireciprocal multiples}

The self-reciprocal feasibility problems have no solution in degrees
\(20\leq D\leq37\). Under the substitution \(x\mapsto-x\), Mossinghoff's
degree-\(38\) multiple of \(\ell^2\)
\cite[Algorithm~3, p.~323]{Mossinghoff2003} becomes the self-reciprocal
Borwein polynomial \(Q_{38}\). Its coefficients are listed in
Table~\ref{tab:l2-self-polynomials}, and its factorization is
\begin{equation}\label{eq:Q38factor}
 Q_{38}(x)=L(x)^2L(-x)\Phi_4(x)\Phi_6(x)\Phi_{12}(x),
\end{equation}
where \(\Phi_n\) denotes the \(n\)-th cyclotomic polynomial.
We therefore obtain
\begin{equation}
 \delta_2^{\sr}(L)=38.
\end{equation}

The antireciprocal feasibility problems have no solution in degrees
\(20\leq D\leq40\). At degree \(41\), an antireciprocal
Borwein multiple exists, whose coefficients are listed in
Table~\ref{tab:l2-anti-polynomials}. Thus
\begin{equation}
 \delta_2^{\ar}(L)=41.
\end{equation}

\subsection{The complete degree spectrum of multiples without symmetry restrictions}

We obtain
\[
 \Dcal_2(L)\cap[34,68]
 =\{34,38\}\cup\{D\in\Z:40\leq D\leq68\}.
\]
For every \(D\geq69\), put \(t=D-34\).  Since \(t>34\), the supports of
\(Q_{34}(x)\) and \(x^tQ_{34}(x)\) are disjoint.  Consequently
\begin{equation}\label{eq:shiftconstruction}
 Q_{34}(x)(1\pm x^t)
 =Q_{34}(x)\pm x^tQ_{34}(x)
\end{equation}
has coefficients in \(\{-1,0,1\}\), degree \(D\), and is divisible by
\(L^2\), for either choice of sign.  Combining this construction with the
classification through degree \(68\) gives
\begin{equation}
 \Dcal_2(L)=\{34,38\}\cup\{D\in\Z:D\geq40\}.
\end{equation}
Thus \(35,36,37\), and \(39\) are the only gaps above the least degree.

\subsection{The degree spectra of self-reciprocal and antireciprocal multiples}

For \(20\leq D\leq76\), the two symmetry-restricted systems give
\begin{align*}
 \Dcal_2^{\sr}(L)\cap[20,76]
   &=\{38\}\cup\{D\in\Z:40\leq D\leq76\},\\
 \Dcal_2^{\ar}(L)\cap[20,76]
   &=\{41\}\cup\{D\in\Z:43\leq D\leq76\}.
\end{align*}
One polynomial for every admitted degree in these two ranges is given in
Tables~\ref{tab:l2-self-polynomials} and \ref{tab:l2-anti-polynomials}.

For every \(D\geq77\), let \(t=D-38\). Since \(Q_{38}\) has degree
\(38\) and \(t>38\), no power of \(x\) occurs in both \(Q_{38}(x)\)
and \(x^tQ_{38}(x)\). Their sum and difference therefore have
coefficients in \(\{-1,0,1\}\). Consequently, the two products
\[
 Q_{38}(x)(1+x^t)
 \qquad\hbox{and}\qquad
 Q_{38}(x)(1-x^t)
\]
are Borwein multiples of \(L^2\) of degree \(D\). The first is
self-reciprocal, while the second is antireciprocal.
It follows that the complete symmetry-restricted spectra are
\begin{equation}\label{eq:square-symmetry-spectra}
 \Dcal_2^{\sr}(L)=\{38\}\cup\{D\in\Z:D\geq40\},
 \qquad
 \Dcal_2^{\ar}(L)=\{41\}\cup\{D\in\Z:D\geq43\}.
\end{equation}

\section{Borwein multiples of \texorpdfstring{$L^3$}{L cubed}}

We now consider \(m=3\). We first tested the unrestricted feasibility
system~\eqref{eq:unrestricted-system} for every degree
\(30\leq D\leq118\). At higher degrees, these unrestricted problems
become computationally demanding, making a systematic degree-by-degree
search impractical. The symmetry constraints of
Section~\ref{subsec:reciprocal-model}, however, reduce the number of
independent variables by approximately half, allowing us to extend the
computations to higher degrees. Since the unrestricted tests exclude
every degree below \(115\), we examined the self-reciprocal and
antireciprocal systems for each \(115\leq D\leq150\). These computations
establish the minimum degrees and the classifications in the indicated
ranges. Further admissible degrees are obtained by shifted constructions.

\subsection{The least degree among multiples without symmetry restrictions}

The unrestricted feasibility problems have no solution for any degree
\[
 30\leq D\leq114.
\]
At degree \(115\), the substitution \(x\mapsto-x\) in Mossinghoff's example
\cite[pp.~322--323]{Mossinghoff2003} gives the self-reciprocal Borwein
polynomial \(Q_{115}\) in Table~\ref{tab:l3-self-polynomials}.
Mossinghoff's factorization~\cite[p.~323]{Mossinghoff2003} becomes
\begin{equation}\label{eq:Q115factor}
\begin{split}
 Q_{115}(x)=L(x)^3&\Phi_2(x)\Phi_3(x)\Phi_4(x)\Phi_{10}(x)
 \Phi_{14}(x)\Phi_{15}(x)\\
 &\times\Phi_{21}(x)\Phi_{22}(x)P_{40}(x),
\end{split}
\end{equation}
where
\begin{align*}
 P_{40}(x)={}&x^{40}-x^{39}-x^{38}+2x^{37}-2x^{35}+x^{34}+x^{33}-x^{32}
 +x^{30}-x^{28}\\
 &+x^{27}-x^{25}+x^{24}-x^{22}+x^{20}-x^{18}+x^{16}-x^{15}
 +x^{13}-x^{12}\\
 &+x^{10}-x^8+x^7+x^6-2x^5+2x^3-x^2-x+1.
\end{align*}
It follows that
\begin{equation}
 \delta_3(L)=115.
\end{equation}

\subsection{The least degrees among self-reciprocal and antireciprocal multiples}

Since \(Q_{115}\) is self-reciprocal and no Borwein multiple exists below
degree \(115\), we obtain
\begin{equation}
 \delta_3^{\sr}(L)=115.
\end{equation}
Thus, imposing self-reciprocity increases the least degree for multiples
of \(L^2\), but leaves it unchanged for multiples of \(L^3\).

At degree \(115\), the antireciprocal feasibility problem has no
solution. At degree \(116\), no Borwein multiple exists, even without
symmetry restrictions. At degree \(117\), both a self-reciprocal and an
antireciprocal multiple exist; their
coefficients appear in Tables~\ref{tab:l3-self-polynomials} and
\ref{tab:l3-anti-polynomials}.  The antireciprocal example is exactly the
sign transform of the degree-\(117\) polynomial exhibited by
Mossinghoff~\cite[p.~323]{Mossinghoff2003}. Since no antireciprocal
multiple exists at degrees \(115\) and \(116\), we obtain
\begin{equation}
 \delta_3^{\ar}(L)=117.
\end{equation}
At degree \(118\), an antireciprocal Borwein multiple exists
(see Table~\ref{tab:l3-anti-polynomials}), whereas no self-reciprocal
multiple exists. Hence
\[
 118\in\Dcal_3^{\ar}(L)\setminus\Dcal_3^{\sr}(L).
\]
Together with the absence of multiples below degree \(115\), this gives
\[
 \Dcal_3(L)\cap[30,118]=\{115,117,118\}.
\]

\subsection{The degree spectra of self-reciprocal and antireciprocal multiples up to degree 150}
\label{subsec:cube-spectrum}

For \(115\leq D\leq150\), self-reciprocal Borwein multiples exist exactly at
the degrees
\begin{equation}\label{eq:cube-selfrec-spectrum-115-150}
\begin{split}
 \Dcal_3^{\sr}(L)\cap[115,150]
 =\{&115,117,120,122,126,127,128,129,\\
    &131,132,133,135,137,138,139,140,141,150\}.
\end{split}
\end{equation}
Antireciprocal Borwein multiples exist exactly at the degrees
\begin{equation}\label{eq:cube-antirec-spectrum-115-150}
\begin{split}
 \Dcal_3^{\ar}(L)\cap[115,150]
 =\{&117,118,123,124,125,127,130,132,133,\\
    &134,135,136,137,138,139,145,146\}.
\end{split}
\end{equation}
The self-reciprocal and antireciprocal spectra overlap in the degrees
\[
 \{117,127,132,133,135,137,138,139\};
\]
at each of these degrees there are distinct self-reciprocal and
antireciprocal multiples.  Tables~\ref{tab:l3-self-polynomials} and
\ref{tab:l3-anti-polynomials} give one polynomial in every admitted degree
of the respective class.

\subsection{Further admissible degrees from shifted constructions}
\label{subsec:cube-shifts}

The polynomials in Tables~\ref{tab:l3-self-polynomials} and
\ref{tab:l3-anti-polynomials} also produce further degrees.  Write
\(Q(x)=\displaystyle\sum_{j=0}^{d}q_jx^j\), with \(q_j=0\) outside
\(0\leq j\leq d\).  The coefficient of \(x^k\) in
\(Q(x)(1\pm x^t)\) is \(q_k\pm q_{k-t}\).  If \(t>d\), the supports
\([0,d]\) and \([t,t+d]\) are disjoint, so both products are automatically
Borwein polynomials.  Taking \(Q=Q_{115}\) and \(t=D-115\) gives a
self-reciprocal product with \(1+x^t\) and an antireciprocal product with
\(1-x^t\) for every \(D\geq231\).

If \(t\leq d\), the supports overlap from degree \(t\) through degree \(d\).
The product is then a Borwein polynomial precisely when none of the integers
\(q_k\pm q_{k-t}\), \(t\leq k\leq d\), equals \(2\) or \(-2\).  Write
\(Q_d^{\sr}\) and \(Q_d^{\ar}\) for the polynomials of degree \(d\) in
Tables~\ref{tab:l3-self-polynomials} and
\ref{tab:l3-anti-polynomials}.  This coefficientwise test gives the
additional degrees

\begin{equation}\label{eq:cube-shifted-degrees}
\begin{split}
 &\{187,188,190,192,195,197,199,200,203,204,205,206,\\
 &\qquad 207,208,209,212,213\}\cup\{215,\ldots,230\}.
\end{split}
\end{equation}
One product for every degree in \eqref{eq:cube-shifted-degrees} is listed
below.  Each product has coefficients in \(\{-1,0,1\}\); degrees at least
\(231\) are already covered by the disjoint-support construction.

\begin{center}
\small
\renewcommand{\arraystretch}{1.12}
\begin{tabular}{@{}r@{\quad}l@{\qquad}r@{\quad}l@{}}
\hline
$D$ & product & $D$ & product \\
\hline
187 & $Q_{115}^{\sr}(x)\bigl(1+x^{72}\bigr)$ & 215 & $Q_{118}^{\ar}(x)\bigl(1+x^{97}\bigr)$ \\
188 & $Q_{115}^{\sr}(x)\bigl(1-x^{73}\bigr)$ & 216 & $Q_{124}^{\ar}(x)\bigl(1+x^{92}\bigr)$ \\
190 & $Q_{115}^{\sr}(x)\bigl(1+x^{75}\bigr)$ & 217 & $Q_{117}^{\ar}(x)\bigl(1-x^{100}\bigr)$ \\
192 & $Q_{115}^{\sr}(x)\bigl(1-x^{77}\bigr)$ & 218 & $Q_{115}^{\sr}(x)\bigl(1-x^{103}\bigr)$ \\
195 & $Q_{115}^{\sr}(x)\bigl(1+x^{80}\bigr)$ & 219 & $Q_{118}^{\ar}(x)\bigl(1+x^{101}\bigr)$ \\
197 & $Q_{117}^{\ar}(x)\bigl(1-x^{80}\bigr)$ & 220 & $Q_{115}^{\sr}(x)\bigl(1+x^{105}\bigr)$ \\
199 & $Q_{118}^{\ar}(x)\bigl(1-x^{81}\bigr)$ & 221 & $Q_{115}^{\sr}(x)\bigl(1+x^{106}\bigr)$ \\
200 & $Q_{115}^{\sr}(x)\bigl(1-x^{85}\bigr)$ & 222 & $Q_{115}^{\sr}(x)\bigl(1-x^{107}\bigr)$ \\
203 & $Q_{117}^{\ar}(x)\bigl(1+x^{86}\bigr)$ & 223 & $Q_{115}^{\sr}(x)\bigl(1-x^{108}\bigr)$ \\
204 & $Q_{117}^{\sr}(x)\bigl(1-x^{87}\bigr)$ & 224 & $Q_{115}^{\sr}(x)\bigl(1+x^{109}\bigr)$ \\
205 & $Q_{115}^{\sr}(x)\bigl(1+x^{90}\bigr)$ & 225 & $Q_{115}^{\sr}(x)\bigl(1+x^{110}\bigr)$ \\
206 & $Q_{115}^{\sr}(x)\bigl(1-x^{91}\bigr)$ & 226 & $Q_{115}^{\sr}(x)\bigl(1-x^{111}\bigr)$ \\
207 & $Q_{122}^{\sr}(x)\bigl(1-x^{85}\bigr)$ & 227 & $Q_{115}^{\sr}(x)\bigl(1-x^{112}\bigr)$ \\
208 & $Q_{117}^{\ar}(x)\bigl(1+x^{91}\bigr)$ & 228 & $Q_{117}^{\sr}(x)\bigl(1-x^{111}\bigr)$ \\
209 & $Q_{117}^{\ar}(x)\bigl(1+x^{92}\bigr)$ & 229 & $Q_{115}^{\sr}(x)\bigl(1+x^{114}\bigr)$ \\
212 & $Q_{117}^{\ar}(x)\bigl(1+x^{95}\bigr)$ & 230 & $Q_{115}^{\sr}(x)\bigl(1-x^{115}\bigr)$ \\
213 & $Q_{115}^{\sr}(x)\bigl(1+x^{98}\bigr)$ & & \\
\hline
\end{tabular}
\end{center}

The overlapping products cover every degree from \(215\) through \(230\);
the disjoint-support construction covers every degree from \(231\) onward.
Consequently
\[
 \{D\in\Z:D\geq215\}\subseteq\Dcal_3(L).
\]
The factor \(1+x^t\) is self-reciprocal and preserves the symmetry type of
\(Q\), whereas \(1-x^t\) is antireciprocal and reverses it.  Below degree
\(231\), the new self-reciprocal degrees obtained in this way are
\[
\begin{split}
 \{&187,190,195,197,199,205,206,209,213,215,216,\\
   &217,218,220,221,223,224,225,227,229,230\},
\end{split}
\]
whereas the new antireciprocal degrees are
\[
\begin{split}
 \{&188,192,195,200,203,204,206,207,208,209,212,213,\\
   &215,216,218,219,220,221,222,223,224,225,226,227,\\
   &228,229,230\}.
\end{split}
\]
The unrestricted degree spectrum of \(L^3\) is therefore reduced to
deciding the following \(55\) degrees:
\[
\begin{split}
 \{&119,121,142,143,144,147,148,149\}\cup\{151,\ldots,186\}\\
 &\cup\{189,191,193,194,196,198,201,202,210,211,214\}.
\end{split}
\]

\section{Conclusion}

The least degrees of Borwein multiples of \(L^2\) and \(L^3\) are now
determined, both without symmetry restrictions and in the
self-reciprocal and antireciprocal classes. The effect of symmetry is
not uniform: imposing self-reciprocity increases the least degree for
multiples of \(L^2\), but leaves it unchanged for multiples of \(L^3\).

Beyond these minima, the degree spectra of \(L^2\) are completely
determined. For \(L^3\), the combination of integer-linear programming
and explicit constructions provides exact classifications in the
ranges considered and establishes the existence of Borwein multiples
in every degree \(D\geq215\), with both symmetry types occurring in
every degree \(D\geq231\).

The remaining questions for \(L^3\) therefore concern only finitely
many degrees. Completing these spectra and extending the investigation
to \(L^4\) are natural directions for further work.

\section*{Appendix: Explicit coefficient vectors}
\phantomsection
\addcontentsline{toc}{section}{Appendix: Explicit coefficient vectors}

All vectors in the tables are ordered by ascending powers of \(x\). Thus a
row \(q_0\ q_1\ \cdots\) begins the coefficient vector of
\(\displaystyle Q(x)=\sum_{j=0}^{D}q_jx^j\). For \(Q_{34}=L^2H_{14}\), the full vector
is displayed. In the self-reciprocal and antireciprocal tables, only the
coefficients \(q_0,\ldots,q_{\lfloor D/2\rfloor}\) are listed. The remaining
coefficients are determined by \(q_{D-j}=q_j\) in the self-reciprocal case
and \(q_{D-j}=-q_j\) in the antireciprocal case. In the latter case, the
central coefficient is necessarily zero when \(D\) is even.

\begingroup
\hbadness=10000
\setlength{\LTleft}{0pt}
\setlength{\LTright}{0pt}
\setlength{\LTpre}{1pt}
\setlength{\LTpost}{1pt}
\setlength{\tabcolsep}{2pt}
\renewcommand{\arraystretch}{0.72}
\begin{longtable}{@{}>{\centering\arraybackslash\fontsize{6}{6.4}\selectfont}p{0.65cm}>{\fontsize{5.2}{5.6}\selectfont\ttfamily\raggedright\arraybackslash}p{14.5cm}@{}}
\caption{The degree-$34$ Borwein multiple $Q_{34}$ of $L^2$.}\label{tab:l2-unrestricted-polynomials}\\
\hline
$D$ & coefficients $q_0,\ldots,q_D$ \\
\hline
\endfirsthead
\endhead
34 & 1 0 -1 -1 0 0 -1 0 1 1 1 0 1 1 0 -1 -1 -1 -1 1 0 -1 1 -1 -1 0 1 1 1 0 0 0 -1 -1 1 \\
\hline
\end{longtable}
\begin{longtable}{@{}>{\centering\arraybackslash\fontsize{6}{6.4}\selectfont}p{0.65cm}>{\fontsize{5.2}{5.6}\selectfont\ttfamily\raggedright\arraybackslash}p{14.5cm}@{}}
\caption{Self-reciprocal Borwein multiples of $L^2$ in every admitted degree from $20$ through $76$.}\label{tab:l2-self-polynomials}\\
\hline
$D$ & listed coefficients $q_0,\ldots,q_{\lfloor D/2\rfloor}$ \\
\hline
\endfirsthead
\multicolumn{2}{c}{\small Table~\thetable\ (continued)}\\
\hline
$D$ & listed coefficients $q_0,\ldots,q_{\lfloor D/2\rfloor}$ \\
\hline
\endhead
\hline
\endfoot
\hline
\endlastfoot
38 & 1 0 -1 0 0 -1 -1 0 1 1 -1 0 1 0 -1 1 0 1 0 0 \\
40 & 1 -1 0 0 0 -1 0 0 1 0 -1 0 1 -1 1 1 1 1 -1 -1 -1 \\
41 & 1 0 0 -1 -1 0 0 0 0 0 0 -1 0 -1 1 1 0 1 0 1 1 \\
42 & 1 0 0 -1 -1 -1 0 0 1 0 0 0 1 0 1 1 1 0 -1 -1 0 0 \\
43 & 1 0 0 -1 0 -1 -1 -1 1 1 1 -1 1 -1 0 -1 1 1 1 1 0 0 \\
44 & 1 -1 -1 0 1 -1 0 1 1 0 0 0 1 0 -1 0 -1 -1 -1 1 -1 1 1 \\
45 & 1 0 -1 -1 0 1 1 -1 -1 0 0 -1 0 0 1 1 -1 -1 1 1 1 1 0 \\
46 & 1 1 -1 -1 -1 0 -1 -1 0 1 0 0 1 1 1 1 0 0 -1 -1 1 0 1 1 \\
47 & 1 0 0 -1 -1 -1 0 1 0 1 0 0 -1 -1 0 1 1 1 0 1 -1 1 0 0 \\
48 & 1 0 -1 -1 0 -1 0 1 1 1 1 -1 0 -1 0 1 0 -1 1 1 -1 0 0 -1 0 \\
49 & 1 -1 0 -1 1 0 0 -1 1 0 0 -1 1 -1 1 -1 1 0 0 0 1 0 -1 0 1 \\
50 & 1 -1 -1 0 0 0 1 0 1 0 0 0 1 -1 0 1 -1 -1 0 0 1 0 -1 0 -1 -1 \\
51 & 1 -1 1 -1 0 -1 0 0 0 0 0 0 0 -1 1 1 0 1 0 1 -1 1 0 0 0 0 \\
52 & 1 0 0 -1 -1 0 -1 0 0 1 0 0 0 0 1 1 1 1 -1 1 0 0 -1 0 1 -1 -1 \\
53 & 1 1 -1 -1 -1 0 -1 -1 0 1 1 -1 0 1 1 1 1 0 0 1 1 -1 0 -1 1 1 -1 \\
54 & 1 -1 1 -1 -1 0 -1 0 0 1 1 1 1 -1 1 1 -1 0 -1 1 -1 0 -1 0 1 -1 -1 1 \\
55 & 1 -1 1 -1 0 -1 0 -1 1 -1 1 1 1 -1 1 -1 0 0 -1 1 1 1 0 -1 -1 1 1 0 \\
56 & 1 0 -1 -1 1 0 -1 -1 1 1 -1 -1 1 1 0 -1 1 1 1 0 0 1 0 0 -1 0 1 -1 -1 \\
57 & 1 0 -1 -1 1 -1 0 -1 1 1 0 -1 1 1 -1 0 0 0 1 1 0 1 -1 0 -1 0 -1 1 0 \\
58 & 1 -1 1 0 -1 -1 0 -1 0 0 0 1 1 -1 1 1 0 0 0 0 -1 0 0 0 1 1 0 1 1 -1 \\
59 & 1 -1 0 0 0 -1 1 -1 1 0 -1 -1 1 0 1 1 0 1 -1 -1 -1 1 0 1 -1 1 -1 0 1 1 \\
60 & 1 -1 0 1 -1 -1 0 0 0 -1 1 1 1 -1 0 1 -1 -1 0 0 1 -1 1 0 1 1 1 0 1 -1 0 \\
61 & 1 -1 0 1 -1 -1 0 -1 1 0 0 1 1 -1 1 0 -1 0 1 -1 -1 -1 1 0 0 1 1 1 1 -1 1 \\
62 & 1 0 -1 -1 1 -1 1 -1 0 0 -1 0 1 0 1 1 1 -1 1 0 1 1 -1 0 -1 -1 0 -1 0 0 1 0 \\
63 & 1 0 -1 0 -1 0 -1 1 1 0 -1 1 0 0 -1 1 1 0 -1 0 1 -1 1 0 -1 0 1 -1 1 0 1 1 \\
64 & 1 -1 1 -1 0 -1 -1 -1 1 1 1 1 1 -1 1 0 -1 0 -1 1 -1 -1 0 1 1 -1 1 0 1 -1 1 -1 0 \\
65 & 1 -1 1 -1 0 -1 -1 0 1 1 0 0 0 -1 1 -1 1 1 0 0 -1 1 1 1 0 -1 0 -1 0 0 1 0 0 \\
66 & 1 -1 -1 0 -1 1 1 1 1 1 -1 -1 0 -1 1 0 -1 -1 -1 -1 0 1 1 0 1 -1 1 -1 1 -1 1 0 -1 0 \\
67 & 1 -1 -1 0 0 -1 1 1 1 1 -1 1 1 -1 0 1 -1 -1 -1 0 0 1 -1 1 -1 -1 0 -1 -1 0 -1 -1 -1 -1 \\
68 & 1 -1 0 -1 1 -1 1 0 1 -1 0 -1 1 -1 1 0 0 -1 -1 1 0 1 1 1 0 0 -1 0 0 0 -1 1 1 0 1 \\
69 & 1 -1 0 -1 0 -1 1 0 1 1 1 0 1 -1 0 0 -1 -1 0 0 -1 1 -1 0 0 1 -1 1 -1 0 0 -1 -1 0 -1 \\
70 & 1 0 -1 0 -1 0 1 -1 1 0 -1 -1 0 0 1 1 -1 1 1 -1 1 -1 1 1 -1 1 0 1 1 -1 0 1 0 1 1 0 \\
71 & 1 -1 0 0 1 -1 -1 -1 0 1 1 0 1 -1 -1 0 -1 1 1 1 0 0 0 -1 1 1 0 1 -1 1 1 -1 1 -1 1 0 \\
72 & 1 -1 1 -1 -1 1 -1 1 -1 0 0 0 0 -1 1 1 -1 1 -1 1 1 1 0 -1 1 0 0 1 0 1 1 -1 0 1 -1 1 0 \\
73 & 1 -1 0 -1 0 1 1 -1 1 -1 1 -1 -1 -1 1 0 -1 -1 1 1 1 -1 1 1 0 0 -1 1 0 -1 1 1 1 1 1 1 1 \\
74 & 1 -1 -1 1 -1 1 1 -1 1 -1 -1 0 0 0 1 0 -1 1 1 0 0 -1 1 -1 -1 1 1 1 1 -1 -1 1 -1 1 0 -1 0 0 \\
75 & 1 0 0 -1 1 -1 -1 -1 0 -1 0 -1 0 1 1 1 1 1 1 1 0 1 0 1 -1 1 1 1 0 1 1 0 0 -1 -1 1 -1 -1 \\
76 & 1 0 -1 -1 1 -1 0 0 -1 1 0 0 1 1 0 1 0 -1 1 1 -1 0 -1 -1 -1 1 1 0 -1 0 1 -1 0 -1 1 0 -1 1 1 \\
\end{longtable}
\begin{longtable}{@{}>{\centering\arraybackslash\fontsize{6}{6.4}\selectfont}p{0.65cm}>{\fontsize{5.2}{5.6}\selectfont\ttfamily\raggedright\arraybackslash}p{14.5cm}@{}}
\caption{Antireciprocal Borwein multiples of $L^2$ in every admitted degree from $20$ through $76$.}\label{tab:l2-anti-polynomials}\\
\hline
$D$ & listed coefficients $q_0,\ldots,q_{\lfloor D/2\rfloor}$ \\
\hline
\endfirsthead
\multicolumn{2}{c}{\small Table~\thetable\ (continued)}\\
\hline
$D$ & listed coefficients $q_0,\ldots,q_{\lfloor D/2\rfloor}$ \\
\hline
\endhead
\hline
\endfoot
\hline
\endlastfoot
41 & 1 -1 0 0 0 -1 0 0 1 0 -1 0 1 -1 0 1 1 1 0 1 0 \\
43 & 1 -1 0 -1 0 0 1 0 1 -1 0 0 1 -1 1 0 0 -1 -1 0 1 0 \\
44 & 1 1 -1 -1 -1 -1 -1 -1 1 1 1 1 1 1 0 1 0 -1 0 0 1 -1 0 \\
45 & 1 -1 0 -1 1 -1 1 0 1 -1 -1 -1 1 0 1 1 1 -1 -1 1 0 1 -1 \\
46 & 1 -1 0 0 -1 0 0 0 1 0 0 1 0 -1 1 1 -1 0 -1 0 0 -1 0 0 \\
47 & 1 -1 0 0 0 -1 1 -1 1 0 -1 -1 1 -1 1 1 0 1 1 0 0 1 0 1 \\
48 & 1 -1 0 0 0 -1 0 0 1 0 0 -1 1 -1 0 0 1 1 1 0 0 0 0 -1 0 \\
49 & 1 -1 -1 1 0 0 -1 0 1 1 -1 -1 1 0 1 0 -1 1 1 -1 -1 -1 1 -1 1 \\
50 & 1 0 0 -1 -1 -1 -1 1 1 0 1 1 0 -1 0 1 0 1 -1 1 0 -1 1 0 -1 0 \\
51 & 1 0 -1 -1 0 0 0 -1 1 1 1 -1 0 1 1 1 -1 -1 0 0 -1 -1 0 1 -1 1 \\
52 & 1 -1 -1 0 1 0 0 0 1 0 -1 -1 1 0 0 0 0 0 -1 0 0 0 1 1 1 1 0 \\
53 & 1 0 -1 0 0 -1 -1 0 1 0 0 0 1 0 0 1 0 0 0 0 1 -1 1 1 0 0 0 \\
54 & 1 0 -1 -1 0 0 0 0 0 1 0 -1 1 1 1 1 -1 -1 -1 0 -1 1 1 0 0 -1 0 0 \\
55 & 1 1 -1 -1 -1 -1 0 -1 0 1 0 1 1 1 0 0 1 0 0 0 1 1 -1 0 -1 1 1 -1 \\
56 & 1 -1 0 0 0 0 -1 0 1 -1 0 -1 1 0 1 0 1 1 0 0 -1 -1 1 -1 1 1 1 1 0 \\
57 & 1 0 0 -1 -1 0 -1 -1 1 1 0 1 0 1 1 -1 -1 1 -1 1 1 0 0 1 -1 0 0 1 1 \\
58 & 1 -1 1 -1 0 -1 0 -1 1 -1 1 0 1 -1 1 1 0 1 0 1 0 0 -1 0 -1 0 -1 1 -1 0 \\
59 & 1 -1 -1 0 -1 1 1 1 1 0 0 0 -1 -1 0 1 -1 1 -1 -1 0 -1 0 -1 0 1 1 -1 0 0 \\
60 & 1 -1 0 -1 1 -1 1 -1 1 0 1 -1 1 -1 1 -1 -1 -1 1 1 -1 1 0 1 0 1 0 0 0 0 0 \\
61 & 1 -1 1 -1 0 -1 0 0 1 -1 -1 1 0 -1 0 1 1 1 -1 1 1 1 -1 1 -1 1 0 1 1 1 1 \\
62 & 1 1 -1 -1 -1 -1 -1 0 1 1 0 0 1 1 -1 1 1 -1 -1 0 1 1 1 1 1 0 -1 0 0 -1 0 0 \\
63 & 1 -1 -1 1 1 -1 -1 0 1 0 -1 0 1 0 0 0 -1 1 1 0 -1 -1 1 1 0 0 1 0 1 0 0 -1 \\
64 & 1 -1 0 0 -1 1 0 0 0 0 -1 -1 0 0 1 1 0 1 0 0 1 0 1 -1 1 -1 0 0 0 -1 0 0 0 \\
65 & 1 -1 0 -1 0 0 1 0 1 0 -1 -1 1 -1 1 1 0 1 -1 -1 1 1 -1 -1 0 -1 0 -1 0 1 1 -1 1 \\
66 & 1 -1 -1 1 1 0 -1 0 0 -1 -1 -1 0 1 1 1 -1 1 1 0 0 -1 1 1 0 1 1 0 0 1 1 -1 1 0 \\
67 & 1 -1 0 1 -1 0 -1 -1 1 -1 0 1 1 0 0 0 -1 1 1 1 1 -1 1 0 0 0 -1 1 0 -1 -1 -1 1 1 \\
68 & 1 -1 1 -1 -1 -1 0 1 0 1 0 1 0 -1 1 1 0 1 -1 0 -1 0 -1 -1 0 -1 1 -1 1 1 0 -1 0 0 0 \\
69 & 1 -1 0 0 -1 1 -1 1 1 -1 -1 -1 0 1 1 0 0 1 0 0 1 0 1 -1 -1 -1 0 0 1 0 1 -1 1 1 -1 \\
70 & 1 0 -1 -1 0 -1 1 0 1 1 0 -1 0 0 1 0 0 -1 1 -1 0 1 0 0 -1 0 1 -1 1 0 -1 0 0 0 0 0 \\
71 & 1 -1 1 -1 -1 0 1 0 -1 -1 1 0 0 -1 1 1 0 1 0 1 0 1 0 0 -1 1 0 0 -1 -1 0 0 1 0 0 1 \\
72 & 1 -1 -1 1 -1 1 0 -1 1 0 0 -1 1 1 1 0 -1 1 1 -1 0 -1 -1 -1 -1 0 -1 1 1 -1 -1 1 -1 0 -1 1 0 \\
73 & 1 -1 0 0 -1 0 0 1 -1 1 -1 1 1 0 1 1 -1 -1 -1 -1 0 1 0 -1 1 0 0 -1 1 0 0 -1 0 1 0 1 -1 \\
74 & 1 -1 0 0 -1 1 -1 1 1 -1 0 -1 0 -1 0 1 0 1 0 1 1 -1 1 -1 1 0 0 0 0 0 1 -1 0 1 1 0 -1 0 \\
75 & 1 -1 -1 -1 0 1 1 1 1 1 1 -1 -1 -1 0 -1 -1 -1 -1 1 0 0 1 0 1 0 0 -1 -1 -1 -1 0 0 0 -1 0 0 -1 \\
76 & 1 0 -1 -1 0 -1 1 -1 1 1 1 1 0 0 -1 1 0 0 -1 -1 -1 0 0 1 0 1 -1 -1 -1 -1 1 1 1 -1 0 -1 1 0 0 \\
\end{longtable}
\begin{longtable}{@{}>{\centering\arraybackslash\fontsize{6}{6.4}\selectfont}p{0.65cm}>{\fontsize{4.8}{5.2}\selectfont\ttfamily\raggedright\arraybackslash}p{14.5cm}@{}}
\caption{Self-reciprocal Borwein multiples of $L^3$ in every admitted degree from $115$ through $150$.}\label{tab:l3-self-polynomials}\\
\hline
$D$ & listed coefficients $q_0,\ldots,q_{\lfloor D/2\rfloor}$ \\
\hline
\endfirsthead
\multicolumn{2}{c}{\small Table~\thetable\ (continued)}\\
\hline
$D$ & listed coefficients $q_0,\ldots,q_{\lfloor D/2\rfloor}$ \\
\hline
\endhead
\hline
\endfoot
\hline
\endlastfoot
115 & 1 -1 -1 0 1 0 0 0 1 0 -1 -1 1 -1 0 1 0 0 1 0 0 -1 0 0 1 0 1 0 -1 0 1 -1 0 0 0 -1 -1 -1 0 0 -1 0 0 0 -1 -1 -1 -1 0 -1 1 0 0 1 0 0 -1 1 \\
117 & 1 -1 0 -1 1 -1 1 0 0 0 -1 -1 1 0 1 1 1 0 1 0 -1 1 -1 -1 0 0 0 0 0 -1 0 -1 0 0 -1 -1 -1 0 -1 1 -1 0 0 -1 -1 0 0 0 0 -1 1 0 0 0 1 -1 -1 -1 -1 \\
120 & 1 -1 0 -1 0 0 1 0 1 0 -1 -1 1 -1 1 1 -1 0 0 1 1 1 -1 1 -1 -1 -1 1 -1 -1 -1 -1 1 0 0 0 -1 0 -1 -1 1 0 1 -1 0 1 1 1 -1 0 0 1 -1 -1 1 1 -1 0 -1 0 0 \\
122 & 1 -1 -1 0 1 0 0 0 0 1 -1 -1 1 1 0 0 -1 0 0 1 -1 0 -1 0 0 1 -1 0 1 0 -1 1 1 1 0 0 -1 0 -1 1 1 -1 -1 -1 0 -1 -1 0 0 0 0 -1 0 0 0 0 -1 -1 0 0 -1 0 \\
126 & 1 -1 0 0 0 -1 0 -1 1 0 0 1 1 -1 1 1 1 1 0 -1 0 -1 -1 -1 -1 0 -1 -1 0 0 1 0 0 0 -1 -1 1 0 -1 1 0 1 1 0 0 -1 0 0 0 -1 0 -1 -1 0 -1 -1 -1 -1 1 0 -1 1 0 0 \\
127 & 1 -1 0 -1 0 0 0 1 1 0 0 0 -1 -1 1 0 1 1 -1 1 -1 0 -1 -1 -1 1 0 -1 0 0 0 0 -1 1 1 1 1 0 0 1 0 -1 0 0 -1 -1 1 -1 0 -1 -1 -1 -1 0 1 1 1 0 1 0 1 0 1 1 \\
128 & 1 0 -1 -1 -1 1 0 0 0 1 1 0 0 0 0 1 -1 1 -1 1 0 0 -1 -1 -1 0 -1 1 -1 1 0 1 0 0 -1 0 0 -1 -1 -1 -1 0 1 0 -1 -1 0 -1 0 0 0 1 0 0 0 1 -1 -1 0 -1 1 1 -1 1 0 \\
129 & 1 -1 0 0 -1 -1 0 1 1 1 1 0 0 -1 0 1 -1 0 0 0 -1 -1 0 -1 0 -1 -1 0 1 0 1 -1 1 0 1 0 1 0 0 1 1 0 0 1 1 1 -1 -1 0 1 -1 1 1 1 0 1 -1 -1 0 1 1 0 -1 -1 -1 \\
131 & 1 -1 0 -1 0 0 1 0 1 0 -1 -1 1 -1 1 0 1 1 -1 -1 1 1 -1 -1 0 1 1 -1 -1 0 -1 0 0 0 0 0 0 0 -1 0 0 1 -1 -1 -1 0 0 -1 0 -1 -1 1 0 0 -1 0 1 0 0 -1 1 1 -1 -1 -1 0 \\
132 & 1 -1 -1 0 0 0 1 1 0 1 0 -1 0 0 -1 0 -1 -1 0 1 -1 1 0 -1 0 1 0 0 1 -1 0 1 0 1 0 1 -1 0 0 1 0 1 0 1 0 1 -1 1 -1 0 -1 0 0 1 1 -1 0 0 0 1 0 1 0 1 0 -1 0 \\
133 & 1 -1 0 -1 0 0 1 -1 1 1 1 -1 0 -1 1 -1 -1 0 1 0 -1 1 0 0 0 0 0 1 0 0 1 0 -1 0 0 0 0 -1 0 -1 -1 0 0 0 -1 -1 -1 -1 0 1 -1 -1 -1 0 -1 -1 0 0 0 0 1 0 1 1 0 1 0 \\
135 & 1 -1 0 -1 1 -1 0 0 1 0 1 -1 1 -1 1 -1 1 0 1 1 -1 -1 0 -1 -1 0 0 -1 0 0 0 0 0 -1 0 1 0 1 0 0 1 -1 0 1 1 0 -1 -1 -1 1 1 1 0 -1 -1 -1 -1 0 1 -1 0 -1 -1 -1 -1 1 -1 -1 \\
137 & 1 -1 0 0 0 -1 0 -1 0 1 1 1 1 -1 0 1 -1 -1 1 1 1 0 -1 0 0 -1 -1 -1 1 0 -1 -1 -1 -1 -1 0 -1 1 -1 -1 0 1 0 1 -1 1 1 0 1 0 1 0 0 0 1 0 -1 -1 -1 0 0 1 0 1 1 0 -1 1 1 \\
138 & 1 0 -1 -1 0 -1 0 1 1 1 0 -1 1 0 -1 1 1 -1 0 1 1 0 -1 -1 -1 -1 0 1 1 -1 -1 -1 0 -1 0 1 -1 -1 -1 -1 1 0 -1 0 -1 1 0 0 1 0 -1 -1 1 1 -1 1 -1 0 0 -1 1 1 0 -1 0 1 -1 1 1 0 \\
139 & 1 -1 0 -1 0 0 1 -1 1 0 1 0 1 -1 1 -1 0 0 -1 -1 1 -1 -1 -1 1 1 1 -1 0 0 1 0 1 -1 0 -1 1 -1 1 -1 0 0 0 -1 1 1 1 -1 1 0 1 0 0 -1 -1 1 1 1 0 0 0 -1 0 1 0 1 -1 -1 -1 -1 \\
140 & 1 -1 -1 1 0 -1 0 1 1 0 -1 -1 1 -1 0 1 0 1 1 0 0 0 1 -1 -1 -1 0 -1 1 0 0 -1 -1 -1 0 -1 0 1 1 -1 1 -1 1 -1 -1 0 0 0 0 1 0 -1 1 -1 1 1 1 1 0 0 0 0 0 1 0 -1 0 -1 -1 0 0 \\
141 & 1 0 -1 -1 0 -1 1 0 0 1 1 -1 1 0 -1 1 1 0 1 0 -1 0 -1 -1 0 1 -1 -1 0 0 0 1 0 -1 -1 -1 -1 0 0 -1 0 0 0 0 -1 0 -1 -1 0 0 1 -1 -1 0 1 1 -1 -1 -1 -1 0 1 0 1 -1 -1 1 1 -1 -1 -1 \\
150 & 1 -1 0 -1 0 0 1 0 1 -1 0 -1 1 0 1 1 0 0 -1 -1 -1 0 0 0 1 1 0 1 -1 -1 -1 0 1 0 -1 -1 0 0 0 -1 0 0 -1 0 0 -1 1 1 1 1 1 0 1 -1 -1 -1 1 -1 -1 0 -1 0 0 1 0 0 -1 -1 1 1 0 1 0 1 0 0 \\
\end{longtable}
\begin{longtable}{@{}>{\centering\arraybackslash\fontsize{6}{6.4}\selectfont}p{0.65cm}>{\fontsize{4.8}{5.2}\selectfont\ttfamily\raggedright\arraybackslash}p{14.5cm}@{}}
\caption{Antireciprocal Borwein multiples of $L^3$ in every admitted degree from $115$ through $150$.}\label{tab:l3-anti-polynomials}\\
\hline
$D$ & listed coefficients $q_0,\ldots,q_{\lfloor D/2\rfloor}$ \\
\hline
\endfirsthead
\multicolumn{2}{c}{\small Table~\thetable\ (continued)}\\
\hline
$D$ & listed coefficients $q_0,\ldots,q_{\lfloor D/2\rfloor}$ \\
\hline
\endhead
\hline
\endfoot
\hline
\endlastfoot
117 & 1 0 -1 0 -1 -1 0 0 0 1 1 1 1 0 0 1 0 0 -1 0 0 0 0 -1 -1 0 0 -1 -1 -1 -1 1 0 0 -1 0 1 -1 0 -1 -1 1 1 -1 -1 0 0 1 0 -1 0 -1 1 -1 -1 -1 -1 1 1 \\
118 & 1 -1 0 0 -1 0 0 0 0 0 1 1 1 -1 1 1 -1 0 -1 0 -1 0 -1 -1 1 0 0 0 -1 0 0 1 0 1 -1 0 0 -1 0 0 0 0 0 1 -1 -1 1 0 0 1 0 1 1 -1 -1 -1 0 -1 0 0 \\
123 & 1 -1 0 0 0 -1 -1 0 1 1 0 1 1 -1 0 0 -1 1 1 1 -1 0 0 0 -1 -1 -1 -1 -1 0 1 1 0 0 -1 -1 -1 0 -1 1 1 -1 -1 0 1 -1 0 0 0 1 1 0 1 0 0 0 1 0 0 0 0 1 \\
124 & 1 -1 -1 0 1 0 0 0 1 0 -1 -1 1 -1 0 0 1 1 0 0 1 -1 -1 0 1 0 1 -1 0 0 -1 -1 1 -1 -1 0 1 0 0 0 0 1 0 -1 0 0 -1 -1 -1 -1 -1 0 -1 -1 -1 0 0 -1 0 0 0 0 0 \\
125 & 1 -1 0 0 -1 0 0 0 0 1 0 0 1 0 1 1 -1 0 -1 0 -1 0 0 -1 1 0 0 -1 1 -1 -1 0 -1 0 0 0 0 -1 1 0 1 1 0 0 0 0 -1 0 1 -1 1 0 0 1 0 -1 0 1 -1 -1 0 1 -1 \\
127 & 1 -1 0 -1 0 0 0 1 1 0 0 -1 1 0 -1 0 0 1 -1 0 0 0 0 -1 0 1 0 1 -1 -1 0 0 1 0 0 1 -1 0 -1 0 0 0 0 0 0 1 0 -1 0 -1 0 -1 -1 -1 0 1 0 -1 1 0 -1 -1 -1 1 \\
130 & 1 -1 0 0 -1 -1 1 0 1 1 1 -1 0 -1 0 0 0 1 1 0 -1 -1 0 -1 0 1 -1 1 0 -1 -1 1 1 0 0 0 -1 -1 0 1 -1 1 -1 0 0 -1 1 0 -1 0 -1 -1 1 1 -1 0 -1 -1 1 -1 1 1 -1 -1 0 0 \\
132 & 1 -1 0 -1 1 -1 0 0 1 0 1 -1 1 -1 1 0 0 0 0 0 0 1 -1 -1 0 0 -1 0 -1 1 1 -1 -1 -1 0 1 1 -1 -1 1 0 -1 1 1 0 1 0 -1 -1 0 0 0 0 0 -1 0 1 -1 -1 -1 1 0 0 0 1 1 0 \\
133 & 1 0 -1 -1 -1 0 0 1 1 1 1 1 0 -1 -1 0 -1 0 0 1 0 -1 -1 -1 -1 1 1 0 -1 0 0 0 0 1 -1 0 -1 1 1 1 1 1 0 0 0 -1 -1 0 0 0 1 -1 1 -1 0 0 1 0 1 -1 -1 -1 0 0 1 0 0 \\
134 & 1 -1 -1 1 -1 1 0 -1 1 1 0 -1 0 -1 1 1 -1 1 1 -1 -1 0 0 0 1 -1 0 0 -1 -1 1 1 -1 1 -1 1 1 -1 0 -1 0 0 0 0 -1 -1 0 0 -1 1 -1 1 1 -1 -1 -1 1 -1 1 1 0 0 -1 -1 -1 0 -1 0 \\
135 & 1 -1 -1 0 1 0 -1 1 0 0 1 0 1 -1 0 1 -1 -1 -1 1 -1 -1 0 0 1 -1 0 -1 -1 1 1 1 1 1 0 1 1 -1 0 -1 1 0 1 1 0 1 -1 0 0 -1 1 -1 0 0 -1 0 0 0 -1 1 0 1 1 0 0 1 -1 0 \\
136 & 1 -1 0 -1 1 -1 0 0 1 0 1 -1 1 -1 1 0 0 -1 1 1 0 0 -1 -1 0 -1 0 0 -1 1 0 0 0 -1 0 1 -1 -1 -1 1 1 0 -1 0 -1 1 1 -1 0 0 0 1 0 1 1 -1 1 1 -1 0 1 1 1 -1 0 1 0 1 0 \\
137 & 1 -1 -1 1 -1 0 1 0 1 1 -1 -1 0 -1 -1 1 0 1 1 -1 1 0 0 -1 1 1 0 0 0 -1 1 -1 -1 -1 0 0 -1 -1 0 1 -1 1 0 -1 0 0 -1 1 -1 0 -1 0 -1 0 -1 0 0 -1 -1 1 1 1 1 -1 0 0 1 -1 0 \\
138 & 1 -1 0 -1 1 -1 0 0 1 0 1 -1 0 0 1 0 1 0 0 1 -1 -1 -1 -1 0 0 -1 -1 1 0 0 -1 0 0 1 0 1 -1 1 0 1 -1 1 0 0 0 0 -1 1 0 1 -1 0 0 1 -1 -1 0 -1 1 -1 0 1 1 -1 0 0 -1 1 0 \\
139 & 1 -1 0 0 -1 0 0 0 1 0 0 0 1 -1 0 1 0 1 -1 0 1 0 0 -1 0 0 -1 0 -1 1 0 0 0 -1 -1 0 -1 -1 -1 0 -1 1 -1 -1 0 -1 0 1 1 1 0 0 -1 0 1 1 1 0 1 -1 0 1 1 1 0 0 1 1 1 -1 \\
145 & 1 -1 0 -1 1 -1 0 0 1 0 1 0 0 -1 1 -1 0 0 0 1 0 0 0 0 0 0 0 -1 -1 1 0 0 1 0 -1 -1 0 -1 -1 0 -1 0 -1 -1 0 0 1 -1 0 -1 0 0 1 0 0 0 1 1 1 0 1 0 1 -1 -1 0 -1 0 -1 0 1 1 1 \\
146 & 1 -1 0 0 -1 0 0 1 -1 1 -1 1 1 -1 1 1 -1 1 -1 1 -1 1 -1 -1 1 -1 1 -1 0 -1 -1 1 0 1 -1 -1 0 -1 0 0 0 0 0 -1 -1 1 -1 -1 -1 -1 1 -1 0 1 1 -1 -1 -1 -1 0 1 -1 1 1 1 0 -1 -1 0 0 1 0 0 0 \\
\end{longtable}
\endgroup

\end{document}